\documentclass[11pt,leqno]{article} 
\usepackage{graphics}
\newtheorem{thm}{Theorem}[section]
\newtheorem{lma}{Lemma}[section]

\newcommand{\beqa}{\begin{eqnarray}}
\newcommand{\eeqa}{\end{eqnarray}}

\newcommand{\epf}{ \hfill$\diamondsuit$ \medskip}

\newcommand{\ds}{\displaystyle}
\newcommand{\beq}{\begin{equation}}
\newcommand{\eeq}{\end{equation}}
\newcommand{\lbl}{\label}
\newcommand{\s}{\; \;}

\newcommand{\ra}{\rightarrow}

\newcommand{\p}{\varphi}

\title{Harmonic oscillators at resonance, perturbed by a non-linear friction force}

\author{
Philip Korman   \\ 
Department of Mathematical Sciences \\ 
University of Cincinnati \\ 
Cincinnati Ohio 45221-0025 \\
and \\
\\
 Yi Li \\
Department of Mathematics  and Statistics\\
Wright State University \\
Dayton OH  45435
}

\date{}

\begin{document}

\maketitle
\begin{abstract} 
This note is an   addendum to the results of A.C. Lazer and P.O. Frederickson \cite{FL}, and A.C. Lazer \cite{L2} on periodic oscillations, with linear part at resonance. We show that a small modification of the argument in \cite{L2} provides a more general result. It turns out that things are different for the corresponding Dirichlet boundary value problem.
 \end{abstract}

\begin{flushleft}
\end{flushleft}

\begin{flushleft}
\end{flushleft}

\section{Introduction}
\setcounter{equation}{0}
\setcounter{thm}{0}
\setcounter{lma}{0}

We are interested in the existence of $2 \pi$ periodic solutions to the problem ($x=x(t)$)
\beq
\lbl{1}
x''+f(x)x'+n^2x=e(t) \,.
\eeq
Here $e(t) \in C(R)$ satisfies $e(t+2\pi)=e(t)$ for all $t$, $f(u) \in C(R)$, $n \geq 1$ is an integer. The linear part, $x''+n^2x=e(t)$, is at resonance, with the null space spanned by $\cos nt$ and $\sin nt$. 
Define $F(x)=\int_0^x f(t) \, dt$. We assume that the finite limits $F(\infty)$ and  $F(-\infty)$ exist, and
\beq
\lbl{2}
F(-\infty)<F(x)<F(\infty)  \s\s \mbox{for all $x$} \,.
\eeq
Define
\[
A_n=\int_0^{2 \pi} e(t) \cos nt \, dt, \s B_n=\int_0^{2 \pi} e(t) \sin nt \, dt \,.
\]

The following theorem was proved in case $n=1$ by A.C. Lazer \cite{L2}, based on P.O. Frederickson and A.C. Lazer \cite{FL}. The paper \cite{FL} was the precursor to the classical works of E.M. Landesman and A.C. Lazer \cite{L}, and A.C. Lazer and D.E. Leach \cite{L}.

\begin{thm}\lbl{thm:1}
The condition 
\beq
\lbl{3}
\sqrt{A_n^2+B_n^2}<2n \left(F(\infty)-F(-\infty) \right)
\eeq
is necessary and sufficient for the existence of $2\pi$ periodic solution of (\ref{1}).
\end{thm}

We provide a proof for all $n$, by modifying the argument in \cite{L2}. 
\medskip

Remarkably, things are different for the corresponding Dirichlet boundary value problem, for which  we derive a necessary condition for the existence of solutions, but show by a numerical computation that  this condition is not  sufficient. Observe that the condition (\ref{3}) depends on $n$, unlike the condition in A.C. Lazer and D.E. Leach \cite{L}.

\section{The proof }
\setcounter{equation}{0}
\setcounter{thm}{0}
\setcounter{lma}{0}

The following elementary lemmas are easy to prove.
\begin{lma}\lbl{lma:1}
Consider a function $\cos (nt-\p)$, with an integer $n$ and any real $\p$. Denote $P=\{t \in (0,2\pi) \, | \, \cos (nt-\p)>0 \}$ and $N=\{t \in (0,2\pi) \,| \, \cos (nt-\p)<0 \}$. Then
\[
\int_P \cos (nt-\p) \, dt=2, \s\s \int_N \cos (nt-\p) \, dt=-2 \,.
\]
\end{lma}

\begin{lma}\lbl{lma:2}
Consider a function $\sin (nt-\p)$, with an integer $n$ and any real $\p$. Denote $P_1=\{t \in (0,2\pi) \, | \, \sin (nt-\p)>0 \}$ and $N_1=\{t \in (0,2\pi) \,| \, \sin (nt-\p)<0 \}$. Then
\[
\int_{P_1} \sin (nt-\p) \, dt=2, \s\s \int_{N_1} \sin (nt-\p) \, dt=-2 \,.
\]
\end{lma}

\noindent
{\bf Proof of the Theorem \ref{thm:1}:} $\;$
1. {\em Necessity}. Given arbitrary numbers $a$ and $b$, we can find a $\delta \in [0,2\pi)$, so that
\[
a \cos nt +b \sin nt=\sqrt{a^2+b^2} \cos (nt-\delta) \,.
\] 
($\cos \delta=\frac{a}{\sqrt{a^2+b^2}}$, $\sin \delta=\frac{b}{\sqrt{a^2+b^2}}$.) We multiply (\ref{1}) by $a \cos nt$, then by $b \sin nt$, integrate and add the results
\beq
\lbl{4}
I \equiv \int_0^{2 \pi} F(x(t))' \cos (nt-\delta) \, dt=\frac{aA_n+bB_n}{\sqrt{a^2+b^2}} \,.
\eeq
Using that  $x(t)$ is a $2\pi$ periodic solution,  and Lemma \ref{lma:2}, we have
\[
I=n\int_0^{2 \pi} F(x(t)) \sin (nt-\delta) \, dt=n\int_{P_1} + n\int_{N_1} <2n \left(F(\infty)-F(-\infty) \right) \,.
\]
Similarly,
\[
I>-2n \left(F(\infty)-F(-\infty) \right) \,,
\]
and so
\[
|I|<2n \left(F(\infty)-F(-\infty) \right) \,.
\]
On the right in (\ref{4}) we have the scalar product of the vector $(A_n,B_n)$ and an arbitrary unit vector. The condition (\ref{3}) follows.
\medskip

\noindent
2. {\em Sufficiency}. We write our equation $\left(x'+F(x) \right)'+n^2x=e(t)$ in the system form

\beqa
\lbl{6}
& x'=-F(x)+y \\
& y'=-n^2x +e(t) \,.\nonumber
\eeqa
Setting $x =\frac1n X$, $y=Y$, we get

\beqa
\lbl{7}
& X'=-nF(\frac1n X)+nY \\
& Y'=-nX +e(t) \,.\nonumber
\eeqa
Let $r(t)=\sqrt{X^2(t)+Y^2(t)}$. Then

\beq
\lbl{8}
r'(t)=\frac{XX'+YY'}{r(t)}=\frac{-nX F(\frac1n X)+e(t)Y}{r(t)}  \,.
\eeq
We see that if $r(t)$ is large, $r'(t)$ is bounded. It follows that there exists $r_0>0$, so that if $|r(0)|>r_0$, then $r(t)>0$ for all $t \in [0,2\pi]$, thus avoiding a singularity in (\ref{8}). Switching to the polar coordinates $X(t)=r(t) \cos \theta (t)$ and $Y(t)=r(t) \sin \theta (t)$, (\ref{8}) becomes

\beq
\lbl{9}
r'(t)=-nF \left( \frac1n r(t) \cos \theta (t) \right) \cos \theta (t)+e(t) \sin \theta (t) \,.
\eeq
We have $\theta (t)=\tan ^{-1} \frac{Y(t)}{X(t)}$, and 
\[
\theta' (t)=\frac{-YX'+XY'}{X^2+Y^2}=\frac{nY F(\frac1n X)-nX^2-nY^2+e(t)X}{X^2+Y^2} \,.
\]
In  polar coordinates
\beq
\lbl{10}
\theta' (t)=-n +\frac{nF \left( \frac1n r(t) \cos \theta (t) \right) \sin \theta (t)}{r(t)}+\frac{e(t)\cos \theta (t)}{r(t)} \,.
\eeq
We denote by $r(t,c,\p)$ and $\theta (t,c,\p)$ the solution of the system (\ref{9}), (\ref{10}) satisfying the initial conditions $r(0,c,\p)=c$ and $\theta (0,c,\p)=\p$.
\medskip

From (\ref{9}) 
\beq
\lbl{11}
r(t,c,\p)=c+O(1), \s \mbox{as $c \ra \infty$}
\eeq
uniformly in $t, \p \in [0,2\pi]$. Then from (\ref{10})
\beq
\lbl{12}
\theta(t,c,\p)=-nt+\p +o(1), \s \mbox{as $c \ra \infty$}
\eeq
uniformly in $t, \p \in [0,2\pi]$. Integrating (\ref{9})
\[
r(2\pi,c,\p)-r(0,c,\p)=\int_0^{2\pi} \left[-nF \left( \frac1n r(t) \cos \theta (t) \right) \cos \theta (t)+e(t) \sin \theta (t) \right] \, dt \,.
\]
We have $\cos \theta (t)=\cos (nt-\p)+o(1)$, and $\sin \theta (t)=\sin (-nt+\p)+o(1)$, as $c \ra \infty$.
Then, in view of (\ref{11})  and Lemma \ref{lma:1}, the integral on the right gets arbitrarily close to
\[
-2n \left(F(\infty)-F(-\infty) \right)+\int_0^{2\pi} e(t) \sin (-nt+\p) \, dt \,,
\] 
for $c$ sufficiently large. Since
\[
\int_0^{2\pi} e(t) \sin (-nt+\p) \, dt =A_n \sin \p-B_n \cos \p<\sqrt{A_n^2+B_n^2} \,,
\]
it follows by our condition (\ref{3}) that
\[
r(2\pi,c,\p)<r(0,c,\p)=c \,,
\]
for $c$ sufficiently large, uniformly in $ \p \in [0,2\pi]$, say for $c>c_1$. Denote $c_2=\max _{c \leq c_1, \, \p \in [0,2\pi]}r(2\pi,c,\p)$, and $c_3=\max(c_1,c_2)$. (Here $r(2\pi,c,\p)$ is computed by using (\ref{6}).) Then $r(2\pi,c,\p) \leq c_3$, provided that $c \leq c_3$. The map $(c,\p) \ra \left(r(2\pi,c,\p),\theta (2\pi,c,\p) \right)$ is   a continuous map of the ball $c \leq c_3$ into itself. By Brouwer's fixed point theorem it has a fixed point, giving us a $2\pi $ periodic solution.
\epf

\section{A boundary value problem}
\setcounter{equation}{0}
\setcounter{thm}{0}
\setcounter{lma}{0}

Consider the Dirichlet problem
\beq
\lbl{14}
x''-F(x)'+x=e(t), \s 0<t<\pi, \s x(0)=x(\pi)=0 \,.
\eeq
Assume that $F(x)$ satisfies (\ref{2}), $e(t) \in C[0,\pi]$.
The linear part has a kernel spanned by $\sin t$. Denote $A=\int _0^{\pi} e(t) \sin t \, dt$. Then from (\ref{14})
\beqa \nonumber
& A=\int _0^{\pi} F(x(t)) \cos t \, dt<F(\infty) \int _0^{\pi/2}  \cos t \, dt+F(-\infty) \int _{\pi/2}^{\pi}  \cos t \, dt \\ \nonumber
& =F(\infty)-F(-\infty) \,. \nonumber
\eeqa
Similarly,
\[
A>F(-\infty)-F(\infty) \,.
\]
We conclude that

\beq
\lbl{15}
|A|<F(\infty)-F(-\infty)
\eeq
is a necessary condition for the existence of solutions.
\medskip

\noindent
It is natural to ask  if the condition (\ref{15}) is sufficient for the existence of solutions. The following numerical computations indicate that the answer is No.
\medskip

\noindent
{\bf Example} We have solved the problem
\beq
\lbl{16}
x''-F(x)'+x=A \sin t+ \sin 2t, \s 0<t<\pi, \s x(0)=x(\pi)=0 \,,
\eeq
with $F(x)=\frac{x}{\sqrt{x^2+1}}$. Here $F(\pm \infty)=\pm 1$, and so the necessary condition for the existence of solutions is $\ds |A|<2$. 
Writing the solution as $x(t)=\xi \sin t+X(t)$, with $\int _0^{\pi} X(t) \sin t \, dt=0$, for each value of $\xi$ we compute the value of $A$ for which the problem (\ref{16}) has a solution with the first harmonic equal to $\xi$, and that solution $x(t)$, see P. Korman \cite{K} for more details. (I.e., we compute the solution curve $(A,x(t))(\xi)$.) In Figure $1$ we draw the curve $A=A(\xi)$. It suggests that there is an $A_0 \approx -0.3$ so that the problem (\ref{16}) has exactly two solutions for $A \in (A_0,0)$, exactly one solution for $A = (A_0,0)$, and no solutions for all other values of $A$. 
The necessary condition $\ds |A|<2$ is definitely not sufficient!

\begin{figure}
\scalebox{0.9}{\includegraphics{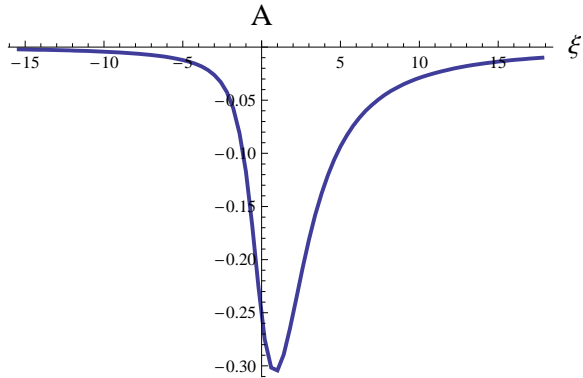}}
\caption{  Solution curve for the problem (\ref{14})}
\end{figure}


\begin{thebibliography}{99}
\bibitem{FL}
P.O. Frederickson and A.C. Lazer,  Necessary and sufficient damping in a second-order oscillator, {\em  J. Differential Equations} {\bf 5}, 262-270 (1969).
\vspace{-0.2cm}

\bibitem{K}
P. Korman,  Global solution curves for boundary value problems, with linear part at resonance, {\em Nonlinear Anal.} {\bf 71}, no. 7-8, 2456-2467  (2009).
\vspace{-0.2cm}

\bibitem{L}
 E.M. Landesman and A.C. Lazer,   Nonlinear perturbations of linear elliptic        
 boundary value problems at resonance, {\em J. Math. Mech.} {\bf 19}, 609-623 (1970).
\vspace{-0.2cm}

\bibitem{L2}
A.C. Lazer, A second look at the first result of Landesman-Lazer type. Proceedings of the Conference on Nonlinear Differential Equations (Coral Gables, FL, 1999), 113-119 (electronic), Electron. J. Differ. Equ. Conf., 5, Southwest Texas State Univ., San Marcos, TX, (2000).
\vspace{-0.2cm}

\bibitem{L1}
A.C. Lazer and D.E. Leach,  Bounded perturbations of forced harmonic oscillators at resonance, {\em Ann. Mat. Pura Appl.}  {\bf 82} (4), 49-68
 (1969). 
\vspace{-0.2cm}

\end{thebibliography}
\end{document}